\documentclass[12pt]{amsart}
\usepackage{graphicx}

\usepackage{amssymb, microtype}

\title[Commentary on Robert Riley's article] {Commentary on Robert Riley's 
article ``A personal account of the discovery of hyperbolic structures on 
some knot complements" }

\begin{document}

\thispagestyle{empty}

\maketitle

\vspace{-30pt}

\[
\hbox{Matthew G.~Brin}^a,\  \hbox{Gareth A.~Jones}^b,\  \hbox{David
Singerman}^c
\]

\noindent \({}^a\)Mathematical Sciences, Binghamton University,
Binghamton, NY, USA.
\newline\noindent
Email: matt@math.binghamton.edu
\newline
\noindent \({}^{b,c}\)Mathematics, University of Southampton,
Southampton, UK
\newline
\noindent Email: \({}^b\)G.A.Jones@soton.ac.uk,
\({}^c\)D.Singerman@soton.ac.uk

\vspace{15pt}

{MSC classification (2010): 
01A60, 01A70} 

\vspace{10pt}


Abstract: We give some background and biographical commentary on the
postumous article \cite{RileyMem} that appears in this journal issue
by Robert Riley on his part of the early history of hyperbolic
structures on some compact 3-manifolds.  A complete list of Riley's
publications appears at the end of this article.

\section{Introduction}

In the mid-1970s the study of the topology of 3-manifolds was
revolutionised by the discovery that many 3-manifolds possess a
hyperbolic structure. This discovery was made, in very different
forms, independently and almost simultaneously, by Robert Riley and
William Thurston, with Riley's results appearing in print first
\cite{MR0425946, MR0412416}.

Riley's approach was algebraic while Thurston's was geometric.
Riley's first results covered a small number of knot complements,
while Thurston's covered large classes of 3-manifolds.  Ultimately a
sweeping conjecture of Thurston \cite{Th2} about the existence of
geometric structures on all 3-manifolds (part of which implies the
Poincar\'e Conjecture) was proven by Perelman.  Riley's earliest
results and conjectures are described in \cite{Th1} as a motivating
factor for Thurston's first result in this area.

Before Riley died in 2000, he wrote a short memoir, 
 describing his recollection of the events leading to his meeting with Thurston in
1976. This was circulated among his colleagues and a few others, but
has lain dormant since his death. A chance conversation between two
people familiar with him from different decades of his life has
revived interest in the memoir, and we have decided to publish it in this 
issue \cite{RileyMem},
along with this commentary which contains some biographical information
and a full bibliography of Riley's
publications.  Riley's wording has not
been altered since, to those who knew him, it reads as pure Riley.
Our only modification to Riley's paper has been to add an abstract,
MSC numbers, keywords, two footnotes and a reference to this paper.



\section{Bob Riley's life and career}

Robert F.~Riley grew up on Long Island in New York State and studied
mathematics at Cornell where he earned his bachelors degree in 1957.
He enrolled in MIT for graduate work with an initial interest in number theory, but was unhappy with the
modern algebraic geometry he was expected to learn there.  Bob spent
some time in industry where he became proficient in the use of
computers.  He regarded himself more as a 19th-century mathematician
with the added advantage of being able to use modern computational
tools.  Much later, Bob proudly showed one of us a letter of rejection he
had received from a reputable British journal saying that they no
longer publish 19th-century mathematics.  

In 1966, Bob moved to Amsterdam.  There he met Brian Griffiths, a
topologist who was Professor of Pure Mathematics at the University of
Southampton.  Brian invited Bob to take a temporary post in
Southampton, which he did in 1968.  In Amsterdam he had become
interested in knot theory and in Southampton he worked on the
representations of knot groups in PSL(2,$\mathbb C$), which is the
group of orientation preserving isometries of hyperbolic 3-space
$\mathbb H^3$.  After some time he realised that, at least for the
figure-eight knot, he was getting a faithful representation, that
the image was a discrete group and that the quotient of $\mathbb
H^3$ by this group was  the figure-eight knot complement.  He had
thus discovered a hyperbolic structure on this knot complement. 
He then showed that the same idea works for several other knots.  
Later, Thurston gave a
necessary and sufficient condition for a knot complement to have a
hyperbolic structure, and wrote that he was 
motivated by Bob's beautiful examples (P. 360 of \cite{Th1}).

Bob discovered these examples with the help of a computer, making
use of his previous industrial experience.  Bob's work was one of
the earliest examples of the extensive use of computers in a branch
of mathematics traditionally dominated by the pure thought and
abstraction method of mathematicians of the first half of the 20th
century.  Note that Bob was working when programs
were submitted to the computer as decks of punched cards.  It should
be mentioned that Thurston, like Riley, was also an outstanding
innovator in computational methods in pure mathematics.

At this time Bob did not have a permanent academic job because he
had left MIT before getting a PhD.  David Singerman agreed to
act as Bob's formal supervisor so that he could be a PhD student at
Southampton.  Bob also obtained funding from the Science Research
Council which financed him for four years.  He obtained a PhD for
his thesis ``Projective representations of link groups", David
Epstein from the University of Warwick being the external examiner.

Bob returned to the USA in 1980, initially joining Thurston at
Boulder and then obtaining a permanent position in the Department of
Mathematical Sciences at Binghamton University in 1982.  He
continued to work there until he died from complications following
heart surgery in March of 2000.

\section{Riley's later mathematics}
In Section 2 of his memoir, Bob wrote the following:

``In December 1991 I used Maple to extend the theorem to algebraic
varieties of $nab$--reps and add some new material.  In 1993 I told
Tomotada Ohtsuki about this, giving no detail, and he promptly found
a better proof and more new material.  I hope to proceed to a joint
paper soon.''

Bob exchanged emails with Ohtsuki, of Kyoto University, in the
mid-1990s, resulting in a short joint manuscript on homomorphisms
between two-bridge knot groups.  This was still incomplete when Bob
died, but it eventually evolved into a joint paper~\cite{MR2484712}
with a third
author, Makoto Sakuma, of Osaka University, who had also
corresponded with Bob (though neither Ohtsuki nor Sakuma had the opportunity to meet
him).

Bob's research sometimes involved quite deep number-theoretic
considerations.  These brought him into contact with the number
theorist Kunrui Yu, now Emeritus Professor at Hong Kong University
of Science and Technology.  In particular, in \cite{MR1041145} Bob
uses the Gel'fond-Baker theory of linear forms in the logarithms of
algebraic numbers to demonstrate the expected growth of the first
homology groups of $k$-sheeted branched covers of $S^3$ branched
over a tame knot.  This paper includes a three-page appendix by Yu.

Yu wrote to us as follows about a visit to Bob in Binghamton in
1990/1.  ``On the second day Bob took me with his very old green
Toyota for a tour.  We visited the factory of Corningware and the
Campus of Cornell University.  The tour was very interesting.  I
found that Bob was a very kind and nice gentleman, and he had very
good sense of humour.''

\subsection{Reminiscences} We remember Bob as an artful eccentric
who practiced his art of bone dry humor, well aware of the effects
he had on his audience.  His accumulated oddities are too numerous
to list and too difficult to explain.  He was fiercely
independent and carved out a life and career that were entirely of
his own making.  He ignored the fashionable, and stuck doggedly to
his own ideas of what was important.  His level of entertainment and
his fellowship were hard to match and he is sorely missed.


\subsection{Thanks}

We would like to thank David Chillingworth, David Epstein, Ross
Geoghegan, Tomotada Ohtsuki, Makoto Sakuma and Kunrui Yu for their
help in producing this commentary.

\bibliographystyle{amsplain}

\def\refname{References}

\bibliographystyle{amsplain}

\def\refname{Publications of Robert Riley}

\end{document}